\documentclass{amsart}
\usepackage{amssymb}
\newtheorem{lemma}{Lemma}[section]
\newtheorem{propos}[lemma]{Proposition}
\newtheorem{theorem}[lemma]{Theorem}
\newtheorem{corol}[lemma]{Corollary}

\newcommand{\CR}{\hbox{{$\mathcal R$}}}
\newcommand{\CE}{\hbox{{$\mathcal E$}}}
\newcommand{\CT}{\hbox{{$\mathcal T$}}}
\newcommand{\CC}{\hbox{{$\mathcal C$}}}
\newcommand{\CS}{\hbox{{$\mathcal S$}}}
\newcommand{\CL}{\hbox{{$\mathcal L$}}}

\newcommand{\R}{\mathbb{R}}
\newcommand{\Z}{\mathbb{Z}}
\newcommand{\C}{\mathbb{C}}

\newcommand{\h}{{\scriptstyle\frac{1}{2}}}
\newcommand{\extd}{{\rm d}}
\newcommand{\del}{\partial}

\newcommand{\image}{{\rm image}}
\newcommand{\ad}{{\rm Ad}}
\newcommand{\isom}{{\cong}}
\newcommand{\eps}{{\epsilon}}
\newcommand{\tens}{\mathop{\otimes}}

\newcommand{\Ad}{{\rm Ad}}
\newcommand{\id}{{\rm id}}
\newcommand{\<}{\langle}
\renewcommand{\>}{\rangle}

\newcommand{\eproof}{$\quad \diamond$\bigskip}
\newcommand{\eqn}[2]{\begin{equation}#2\label{#1}\end{equation}}
\newcommand{\und}[1]{{\underline {#1}}}

\newcommand{\rcross}{{\triangleright\!\!\!<}}
\newcommand{\lcross}{{>\!\!\!\triangleleft}}
\newcommand{\lbiprod}{{>\!\!\!\triangleleft\kern-.33em\cdot}}
\newcommand{\rbiprod}{{\cdot\kern-.33em\triangleright\!\!\!<}}

\begin{document}

\title[NONCOMMUTATIVE DIFFERENTIALS AND
YANG-MILLS ON $S_N$]{\rm \large NONCOMMUTATIVE DIFFERENTIALS AND
YANG-MILLS ON PERMUTATION GROUPS $S_N$}
\author{S. Majid}
\address{School of Mathematical Sciences\\
Queen Mary, University of London\\ 327 Mile End Rd,  London E1
4NS, UK.}

\thanks{This paper is in final form and no version of it will
be submitted for publication elsewhere} \subjclass{58B32, 58B34,
14N15} \keywords{noncommutative geometry, braided categories,
finite groups, quantum groups, flag variety, orbit method}

\date{5/2001; revised 1/2003}%
%\dedicatory{}%
%\commby{}%

\maketitle

\begin{abstract}
\noindent We study noncommutative differential structures on the
group of permutations $S_N$, defined by conjugacy classes. The
2-cycles class defines an exterior algebra $\Lambda_N$ which is a
super analogue of the Fomin-Kirillov algebra $\CE_N$ for Schubert
calculus on the cohomology of the $GL_N$ flag variety.
Noncommutative de Rahm cohomology and moduli of flat connections
are computed for $N<6$. We find that flat connections of
submaximal cardinality form a natural representation associated to
each conjugacy class, often irreducible, and are analogues of the
Dunkl elements in $\CE_N$. We also construct $\Lambda_N$ and
$\CE_N$ as braided groups in the category of $S_N$-crossed
modules, giving a new approach to the latter that makes sense for
all flag varieties.
\end{abstract}

\section{Introduction}

In recent years there has been developed a fully systematic
approach to the noncommutative differential geometry on (possibly
noncommutative) algebras, starting with differential forms on
quantum groups\cite{Wor:dif} and including principal bundles with
Hopf algebra fiber, connections and Riemannian structures, etc,
see \cite{Ma:rief} or our companion paper in the present volume
for a review. These constructions successfully extend conventional
concepts of differential geometry to the $q$-deformed case such as
$q$-spheres and $q$-coordinate rings of quantum groups.

However, this constructive noncommutative geometry can also be
usefully specialised to finite-dimensional Hopf algebras and from
there to finite groups, where differentials and functions
noncommute (even though the functions themselves commute). Indeed,
one has then a rich `Lie theory of finite groups' complete with
differentials, Yang-Mills theory, metrics and Riemannian
structures. If $k(G)$ denotes the functions on the finite group,
then the differential structures are defined by exterior algebras
of the form $\Omega=k(G).\Lambda$ where $\Lambda$ is the algebra
of left-invariant differential forms. These in turn are determined
by conjugacy classes. The case of the symmetric group $S_3$ of
permutations of 3 elements, with its 2-cycle conjugacy class, was
fully studied in \cite{Ma:rief} and \cite{MaRai:ele}. Among other
results, it was shown that $S_3$ has the same noncommutative de
Rham cohomology as the quantum group $SL_q(2)$ (or 3-sphere
$S_q^3$ in a unitary setting). The moduli space of flat $U(1)$
connections on $S_3$ is likewise nontrivial and was computed. The
goal of the present article is to extend some of these results to
higher $S_N$, with some results for all $N$ and others by explicit
computation for $N<6$.

In particular, we make a thorough study of the invariant exterior
algebra $\Lambda=\Lambda_N$ with the 2-cycles calculus, and
explain its close connection with other algebras in mainstream
representation theory (in Schubert calculus) and algebraic
topology. Our first result is a description of the algebra
$\Lambda_N$  as generated by $\{ e_{(ij)}\}$ labelled by 2-cycles
with relations
\[ e_{(ij)}\wedge e_{(ij)}=0,\quad e_{(ij)}\wedge e_{(km)}
+e_{(km)}\wedge e_{(ij)}=0\]
\[e_{(ij)}\wedge e_{(jk)}+e_{(jk)}\wedge e_{(ki)}
 +e_{(ki)}\wedge e_{(ij)}=0\] where
$i,j,k,m$ are distinct. We consider $e_{(ij)}= e_{(ji)}$ since
they are labelled by the same 2-cycle. Our first observation is
that $\Lambda_N$ has identical form to the noncommutative algebra
$\CE_N$ introduced in \cite{FomKir:qua} with generators $[ij]$ and
relations \[ [ij]=-[ji],\quad [ij]^2=0,\quad
[ij][km]=[km][ij],\quad [ij][jk]+[jk][ki]+[ki][ij]=0\] for
distinct $i,j,k,m$. The main difference is that our $e_{(ij)}$ are
symmetric and partially anticommute whereas the $[ij]$ are
antisymmetric and partially commute. We will show that many of the
problems posed in \cite{FomKir:qua} and some of the results there
have a direct noncommutative-geometrical meaning in our super
version. For example, the algebra $\CE_N$ has a subalgebra
isomorphic to the cohomology of the flag variety associated to
$GL_{N}$ and among our analogous results we have a subalgebra of
flat connections with constant coefficients. These results are in
Section~3, with some further metric aspect on Section~5. Moreover,
using our methods we obtain several new results about the algebras
$\CE_N$. These are in Section~6. The first and foremost is our
result that the $\CE_N$ are braided groups or Hopf algebras in
braided categories. We conjecture that as such they are self-dual
and show that this unifies and implies several disparate
conjectures in \cite{FomKir:qua}. We show that the extended
divided-difference operators $\Delta_{ij}$ in that paper are
indeed the natural braided-differential operators on any braided
group, and that the cross product Hopf algebras in
\cite{FomPro:fib} are the natural bosonisations. We also prove
that if $\CE_N$ is finite-dimensional then it has a unique element
of top degree. Our approach works for all flag varieties
associated to other Lie algebras with Weyl groups beyond $S_N$.

The reasons for the close relation between $\Lambda_N$ and $\CE_N$
is not known in detail but can be expected to be something like
this: the flag variety has a cell decomposition labelled by $S_N$
and its differential geometric invariants should correspond in
some sense to the noncommutative discrete geometry of the
`skeleton' of the variety provided by the cell decomposition. One
can also consider this novel phenomenon as an extension of
Schur-Weyl duality. Let us also note the connnection between flag
varieties and the configuration space $C_N(d)$ of ordered
$N$-tuples in $\R^d$ with distinct entries, as emphasized in the
recent works of Lehrer, Atiyah and others. Its cohomology ring in
the Arnold form can be written as generated by $d-1$-forms
$E_{ij}$ labelled by pairs $i\ne j$ in the range $1,\cdots,N$ with
relations \cite{Coh}\cite{Leh}
\[ E_{ij}=(-1)^d E_{ji},\quad E_{ij}E_{km}=(-1)^{d-1}E_{km}E_{ij},
\quad E_{ij}E_{jk}+E_{jk}E_{ki}+E_{ki}E_{ij}=0\] for all $i,j,k,m$
not necessarily distinct (to the extent allowed for the labels to
be valid). We have rewritten the third relation in the required
suggestive form using the first two relations. We see that
$H(C_N(d))$ is precisely a graded-commutative version of the
algebra $\Lambda_N$ if $d$ is even and of $\CE_N$ extended by
dropping the $[ij]^2=0$ relation if $d$ is odd. Thus one can say
that the noncommutative geometry of $S_N$ and the extended
Fomin-Kirillov algebra $\CE_N$ together `quantize' the cohomology
of this configuration space in the sense that some of the
graded-commutativity relations are dropped.

In the preliminary Section~2 we recall the basic ingredients of
the constructive approach to noncommutative geometry (coming out
of quantum groups) that we use. Its relation to other approaches
such as \cite{Con:geo} is only partly understood, see
\cite{Ma:rief}. In Section~4 of the paper we look at other
differential calculi on symmetric groups as defined by other
conjugacy classes. To be concrete we look at $S_4,S_5$ and compute
moduli of flat connections with constant coefficients. The result
suggests the (incomplete) beginnings of an approach to construct
an irreducible representation associated to each conjugacy class
by noncommutative-geometrical means and in a manner that would
make sense in principle for any finite group $G$. Since all of the
geometry is $G$-equivariant there is plenty of scope to associate
representations; here we explore one such method and tabulate the
results.

\section{Preliminaries on noncommutative differentials}

Noncommutative differential geometry works over a general unital
(say) algebra $A$. The main idea is to define the differential
structure by specifying an $A-A$-bimodule $\Omega^1$ of `1-forms'
equipped with an exterior derivative $\extd:A\to\Omega^1$ obeying
the Leibniz rule. When $A$ is a Hopf algebra there is a natural
notion of $\Omega^1$ bicovariant \cite{Wor:dif} and in this case
it can be shown that $\Omega^1=A.\Lambda^1$ (a free left
$A$-module) where $\Lambda^1$ is the space of left-invariant
1-forms. This space has the natural structure of a right
$A$-crossed module (in the case of $A$ finite-dimensional it means
a right module over the right quantum double of $A$) and as a
result a braiding operator $\Psi:\Lambda^1\tens\Lambda^1\to
\Lambda^1\tens\Lambda^1$ obeying the Yang-Baxter equations. This
can be used to define the wedge product between invariant 1-forms
in such a way that they `skew-commute' with respect to $\Psi$. The
naive prescription is a quadratic algebra $\Lambda_{quad}$ but
there is also a more sophisticated Woronowicz prescription
$\Lambda_w$; in both cases the exterior algebra $\Omega$ is
defined as freely generated by these over $A$. Here \eqn{Lambdaw}{
\Lambda_{quad}=T\Lambda^1/\ker(\id-\Psi),\quad
\Lambda_w=T\Lambda^1/\oplus_n\ker A_n} as quotients of the tensor
algebra, where the Woronowicz antisymmetrizer is \eqn{An}{
A_n=\sum_{\sigma\in
S_n}(-1)^{l(\sigma)}\Psi_{i_1}\cdots\Psi_{i_{l(\sigma)}}:(\Lambda^1)^{\tens
n}\to (\Lambda^1)^{\tens n}.} Here $\Psi_i\equiv\Psi_{i,i+1}$
denotes $\Psi$ acting in the $i,i+1$ place and
$\sigma=s_{i_1}\cdots s_{i_{l(\sigma)}}$ is a reduced expression
in terms of simple reflections. There is also an operator
$\extd:\Lambda^1\to \Lambda^2$ which extends to the entire
exterior algebra with $\extd^2=0$ and defines the noncommutative
de Rham cohomology as closed forms modulo exact.
\begin{propos} \cite{Ma:eps}
$A_n=[n;-\Psi]!$ where
\[ [n;-\Psi]=\id-\Psi_{12}+\Psi_{12}\Psi_{23}
+\dots+(-1)^{n-1}\Psi_{12}\cdots\Psi_{n-1,n}\] are the braided
integer matrices and $[n;-\Psi]!=(\id\tens[n-1;-\Psi]!)[n;-\Psi]$.
\end{propos}
This is a practical method to compute the $A_n$, which we will
use. It comes from the author's theory of braided binomials (or
sometimes called braided shuffles) introduced in \cite{Ma:fre}.
See also the later works \cite{Ros}\cite{BesDra:dif}. For example,
\begin{eqnarray*} [3;-\Psi]!\kern -15pt &&=(\id\tens [2;-\Psi])[3;-\Psi]
=(\id-\Psi_{23})(\id-\Psi_{12}+\Psi_{12}\Psi_{23})\\
&&=\id-\Psi_{12}-\Psi_{23}+\Psi_{12}\Psi_{23}+\Psi_{23}\Psi_{12}
-\Psi_{23}\Psi_{12}\Psi_{23}=A_3.\end{eqnarray*}

For other formulae it is enough for our purposes to specialise
directly to finite sets and finite groups.  We work over a field
$k$ of characteristic zero. Let $A=k(\Sigma)$ a finite set. Then
the differential structures are easily seen from the axioms to
correspond to subsets $E\subset \Sigma\times\Sigma-{\rm diag}$ of
`allowed directions'. Thus \eqn{OmegaSig}{ \Omega^1={\rm
span}\{\delta_x\tens\delta_y|\ (x,y)\in E\},\quad \extd
f=\sum_{(x,y)\in E} (f(y)-f(x)) \delta_x\tens\delta_y} where
$\delta_x$ is the Kronecker delta-function. Note that
$\delta_x\tens\delta_y=\delta_x\extd\delta_y$ for all $(x,y)\in
E$. This result for finite sets is common to all approaches to
noncommutative geometry, e.g. in \cite{Con:geo}. If $\Sigma=G$ is
a finite group then a natural choice of $E$ is given by
\eqn{Econj}{ E=\{(x,y)\in G\times G|\ x^{-1}y\in \CC\}} for any
subset $\CC$ not containing the group identity $e$. Such a
calculus is manifestly invariant under translation by $G$ and all
covariant differential calculi are of this form. Bicovariant ones
(as above) are given precisely by those $\CC$ which are
$\Ad$-stable, so that $E$ is invariant from both sides.  The
`simple' such differential structures (with no proper quotient)
are classified precisely by the nontrivial conjugacy classes. They
take the form of a free left $k(G)$-module \eqn{OmegaG}{
\Omega^1=k(G)\cdot{\rm span}\{e_a|\ a\in \CC\},\quad\extd
f=\sum_{a\in \CC}(R_a(f)-f) e_a ,\quad e_a f=R_a(f) e_a} where
$R_a(f)(g)=f(ga)$ denotes right translation and, explicitly,
$e_a=\sum_{g\in G} \delta_g \extd \delta_{ga}$. Such formulae
follow at once from Woronowicz's paper as a special case. An early
study of this case, in the physics literature, is in \cite{Bre}.

Moreover, in the case of a finite group $G$, a right
$k(G)$-crossed module is the same thing (by evaluation) as a left
$G$-crossed module in the sense of Whitehead, i.e. a $G$-graded
$G$-module with the degree map $|\ |$ from the module to $kG$
being equivariant (where $G$ acts on $kG$ by $\ad$),
see\cite{Ma:book}. The particular crossed module structure on
$\Lambda^1=k\CC$ and induced braiding are \eqn{PsiG}{|e_a|=a,\quad
g.e_a=e_{gag^{-1}},\quad \Psi(e_a\tens e_b)=e_{aba^{-1}}\tens
e_a.}

\begin{propos} \cite{Ma:rief}  For each $g\in G$, consider the set
$\CC\cap g \CC^{-1}$. This has an automorphism $\sigma(a)=a^{-1}g$
corresponding to the braiding under the decomposition $k\CC\tens
k\CC=\sum_g k(\CC\cap g\CC^{-1})$. Hence if $V_g=(k\CC\cap g
\CC^{-1})^\sigma$ (the fixed subspace) has basis
$\{\lambda^{(g)\alpha}\}$, the full set of relations of
$\Lambda_{quad}$ are
\[\forall g\in G:\quad \sum_{a,b\in\CC,\, ab=g} \lambda_a^{(g)\alpha}
e_ae_b=0.\] \end{propos}

These are also the relations of $\Omega_{quad}$ over $k(G)$.
Meanwhile, the exterior derivative is provided by \eqn{extdG}{
\extd e_a=\theta e_a+e_a\theta,\quad \theta=\sum_{a\in\CC} e_a.}
It follows that $\extd$ is given in all degrees by
graded-commutation with the 1-form $\theta$. It is easy to see
that it obeys $\theta^2=0$ and $\extd\theta=0$ and that $\theta$
is never exact (so the noncommutative de Rham cohomology $H^1$
always contains the class of $\theta$).

\section{2-Cycle differential structure on $S_N$}

It is straightforward to compute the quadratic exterior algebra
for $G=S_N$ from the above definitions. We are particularly
interested in the invariant differential forms since these
generate the full structure over $k(G)$. In this section, we take
the differential structure defined by the conjugacy class $\CC$
consisting of 2-cycles described as unordered pairs $(ij)$ for
distinct $i,j\in\{1,\cdots,N\}$.

\begin{propos} The quadratic exterior algebra $\Lambda_N\equiv
\Lambda_{quad}(S_N)$ for
the 2-cycles class is the algebra with generators $\{ e_{(ij)}\}$
and relations
\[ (i)\quad e_{(ij)}^2=0,\quad (ii)\quad e_{(ij)}e_{(km)}+e_{(km)}
e_{(ij)}=0\]
\[ (iii)\qquad e_{(ij)}e_{(jk)}+e_{(jk)}e_{(ik)}+e_{(ik)}e_{(ij)}=0\]
where $i,j,k,m$ are distinct.
\end{propos}
\proof There are three kinds of elements $g\in G$ for which
$\CC\cap g\CC^{-1}$ is not empty. These are (i) $g=e$, in which
case $\sigma$ is trivial and $V_e=k\CC$. This gives the relations
(i) stated; (ii) $g=(ij)(km)$ where $i,j,k,m$ are disjoint. In
this case $\CC\cap g\CC^{-1}$ has two elements $(ij)$ and $(km)$,
interchanged by $\sigma$. The basis of $V_{(ij)(km)}$ is
1-dimensional, namely $(ij)+(km)$ and this gives the relation (ii)
stated; (iii) The element $g=(ij)(jk)$ where $i,j,k$ are disjoint.
Here $\CC\cap g\CC^{-1}$ has 3 elements $(ij),(jk),(ik)$
cyclically rotated by $\sigma$. The invariant subspace is
1-dimensional with basis $(ij)+(jk)+(ik)$ giving the relation
(iii). \eproof

We note that
\[ \dim(\Lambda_N^1)=\left({N\atop 2}\right),\quad
\dim(\Lambda_N^2)={N(N-1)(N-2)(3N+7)\over 24}\] which are the same
dimensions as for the algebra $\CE_N$ in \cite{FomKir:qua}. The
first of these is the `cotangent dimension' of the noncommutative
manifold structure on $S_N$. It is more or less clear from the
form of the two algebras that their dimensions coincide in all
degrees (and for $N=3$ they are actually isomorphic). We have
computed these dimensions for the exterior algebra for $N<6$ using
the explicit form of the braiding $\Psi$ defining the algebra, and
they indeed coincide with the corresponding dimensions for $\CE_N$
listed in \cite{FomKir:qua}. These data are listed in Table~1 with
the compact form of the Hilbert series taken from
\cite{FomKir:qua} (for $S_4, S_5$ only the low degrees have been
explicitly verified by us). Also of interest is the top degree $d$
in the last column. From our noncommutative geometry point of view
this is the `volume dimension' of the noncommutative manifold
structure where the top form plays the role of the volume form.
Note that the cotangent dimension and volume dimension need not
coincide even though they would do so in classical geometry. Also
note (thanks to a comment by R. Marsh) that these volume
dimensions are exactly the number of indecomposable modules of the
preprojective algebra of type $SL_N$. The latter is a quotient of
the path algebra of the doubled quiver of the associated oriented
Dynkin diagram and a module means an assignment of `parallel
transport' operators to arrows of the quiver, i.e. some kind of
`connection'. A classic theorem of Lusztig-Kashiwara-Saito states
that there is a 1-1 correspondence between the irreducible
components of its module variety with fixed dimension vector and
the canonical basis elements of the same weight. Since the
representation theory for the next preprojective algebra in the
series is tame but infinite, we therefore expect $\Lambda_6$ and
higher to be infinite-dimensional, and similarly for $\CE_6$ and
higher.

The infinite-dimensionality or not of $\CE_6$ has been posed in
\cite{FomKir:qua}, where it was conjectured that if finite
dimensional then the top form should be unique (we will prove this
in Section~6) and that the Hilbert series should have a symmetric
increasing and decreasing form.  Without proving this second
conjecture here, let us outline a noncommutative-geometric
strategy for its proof. Namely, the Woronowicz quotient
$\Lambda_w$ by its very construction will be nondegenerately
paired with a similar algebra $\Lambda^*_w$ of `skew tensor
fields' (see the Appendix). Moreover, if finite dimensional, and
in the presence of a nondegenerate metric (see in Section~5) we
then expect Hodge * isomorphisms $\Lambda_w^m\to\Lambda_w^{d-m}$
and ultimately an increasing-decreasing symmetric form of the
Hilbert series for $\Lambda_w$ as familiar in differential
geometry. All of this was concretely demonstrated for $S_3$ in
\cite{MaRai:ele}. The main ingredient missing then is that
$\Lambda_N$ is the quadratic quotient whereas the Woronowicz one
could in principle be a quotient of that. The same strategy and
considerations apply to $\CE_N$. As a step we have,

\begin{table}
\[\begin{array}{c|ccccccc}
{\rm dim}& \Omega^0& \Omega^1& \Omega^2 &\Omega^3 &\Omega^4& {\rm
Hilbert\ polynomial}(q)& {\rm Top\ degree}\\ \hline
S_2 & 1& 1 & & & & [2]_q& 1\\
S_3 & 1 & 3 & 4 & 3& 1& [2]_q^2[3]_q& 4\\
S_4 &1 &6&19&42&71&[2]_q^2[3]_q^2[4]_q^2  & 12   \\
S_5 &1 & 10 & 55 & 220 & 711 & [4]_q^4[5]_q^2[6]_q^4& 40\\
\end{array}\]
\caption{Dimensions and Hilbert polynomial for the exterior
algebras $\Lambda_{quad}$ for $N<6$ as for $\CE_N$ in
\cite{FomKir:qua}. Here $[n]_q=(q^n-1)/(q-1)$.}
\end{table}

\begin{theorem} For the 2-cycle differential calculus on $S_N$,
$\Lambda_w=\Lambda_{quad}$ in degree $<4$ (we conjecture this for
all degrees).
\end{theorem}
\proof  We will decompose the space $k\CC\tens k\CC\tens
k\CC=V_3\oplus V_2\oplus V_1\oplus V_0$ where each $V_i$ is stable
under the braiding operators $\Psi_{12},\Psi_{23}$. Since $A_3$
can be factorised either through $\id-\Psi_{12}$ or
$\id-\Psi_{23}$, its kernel contains that of these operators. So
it suffices to show on each $V_i$ that the dimension of the kernel
of $A_3$ equals the dimension of the sum of the kernels of
$\id-\Psi_{12},\id-\Psi_{23}$. We say $a\sim b$ if the 2-cycles
$a,b$ have exactly one entry in common and $a\perp b$ if disjoint.
We then decompose $\CC\times\CC\times\CC$ as follows. For $V_0$ we
take triples $(a,b,c)$ which are pairwise either $\perp$ or equal,
but not all three equal. Here the braiding is trivial. For $V_1$
we take triples where two pairs are mutually $\perp$ and one is
$\sim$. It suffices to let the totally disjoint element be fixed,
say (45) and the others to have entries taken from a fixed set,
say $\{1,2,3\}$ (i.e. $V_1$ is a direct sum of stable subspaces
spanned by basis triples with these properties fixed). On such a
9-dimensional space one may compute $\Psi_{12},\Psi_{23}$
explicitly and verify the required kernel dimensions (for $A_3$ it
is 7). For $V_2$ we take triples $(a,b,c)$ where two pairs are
$\sim$ and one is $\perp$, or where all three pairs are $\sim$
through the same entry occurring in all three 2-cycles. This time
it suffices to take entries from $\{1,2,3,4\}$, say, and verify
the kernels on such a 16-dimensional subspace (for $A_3$ it is 11
dimensional). For $V_3$ we take triples which are pairwise either
$\sim$ or $=$, excluding the special subcase of three $\sim$ used
in $V_2$. Here it suffices to take entries from $\{1,2,3\}$ and
the braidings become as for $S_3$, where the result is known from
\cite{MaRai:ele}. The remaining type of triple, where there is one
pair $\perp$, one $=$ and one $\sim$, is not possible. \eproof

The absence of additional cubic relations strongly suggests that
the Woronowicz exterior algebra on $S_N$ coincides with the
quadratic one in all degrees (this is known for $N=2,3$ by direct
computation). In view of the above theorem, we continue to work
throughout with the quadratic exterior algebra. On the other hand,
it should be stressed that we expect $\Lambda_w=\Lambda_{quad}$ to
be a special feature of $S_N$. The evidence for this is that one
may expect a kind of `Schur-Weyl duality' between the
noncommutative geometry of the finite group on one side and that
of the classical or quantum group on the other. And on the quantum
group side it is known that the Woronowicz exterior algebra of
$SL_q(N)$ coincides with the quadratic one for generic $q$, but
not for the other classical families. Therefore for other than the
$SL_N$ series we would expect to need to work with the
nonquadratic $\Lambda_w$ and likewise propose a corresponding
nonquadratic antisymmetric version generalising the $\CE_N$.

Next, for any differential graded algebra we define cohomology as
usual, namely closed forms modulo exact. It is easy to see that $
H^0(S_N)=k.1$ for all $N$.

\begin{propos} The first noncommutative de Rham cohomology
of $S_N$ at least for $N<6$ with the 2-cycle differential
structure is
\[ H^1(S_N)=k.\theta\]
\end{propos} \proof This is done by direct
 computation of the dimension of the kernel of $\extd$,
along the same lines as in \cite{MaRai:ele}, after which the
result follows. We expect that in fact $H^1(S_N)=k.\theta$ for all
$N$, but the general proof requires some elaboration.\eproof

It follows from Poincar\'e duality that $H^2(S_3)=0$ and
$H^3(S_3)=k$, $H^4(S_3)=k$ as computed explicitly in
\cite{MaRai:ele}, which is the same as for $SL_q(2)$ and gives
some small evidence for the Schur-Weyl duality mentioned above (up
to a shift or mismatch in the rank).  Next, beyond the cohomology
$H^1$ is a nonlinear variant which can be called `$U(1)$
Yang-Mills theory'. Here a connection or gauge field is again a
1-form $\alpha\in \Omega^1$. But rather than modulo exact 1-forms
we are interested in working modulo the gauge transformation
\[ \alpha\mapsto u\alpha u^{-1}+u\extd u^{-1}\]
for invertible $u$ in our coordinate algebra. The covariant
curvature of a connection is
\[ F(\alpha)=\extd \alpha+ \alpha^2\]
and transforms by conjugation. This is like nonAbelian gauge
theory but is nonlinear even for the $U(1)$ case because the
differential calculus is noncommutative.

\begin{propos} For the 2-cycle differential calculus on $S_N$,
\[ \alpha_i=-\theta_i,\quad \theta_i=\sum_{j\ne i} e_{(ij)}\]
are flat connections with constant coefficients. The 1-forms
$\theta_i$ obey $\theta_i\theta_j+\theta_j\theta_i=0$ for $i\ne
j$.
\end{propos}
\proof We first check the anticommutativity for $i\ne j$. In the
sum
\[ \theta_i\theta_j+\theta_j\theta_i=\sum_{k\ne i,l\ne j} e_{(ik)}
e_{(jl)}+e_{(jl)}e_{(ik)}\] only the cases where $i,j,k,l$ are not
distinct contribute due to relation (ii) in Proposition~3.1.
Likewise the terms where $(ik)=(jl)$ do not contribute by (i).
There are three remaining and mutually exclusive cases: $k=l$, or
$k=j$ or $l=i$. Relabelling the summation variable $k$ in each
case we have,
\[\theta_i\theta_j+\theta_j\theta_i=\sum_{k\ne i,j}
e_{(ik)}e_{(jk)}+e_{(jk)}e_{(ik)}+e_{(ij)}e_{(jk)}
+e_{(jk)}e_{(ij)}+e_{(ik)}e_{(ji)}+e_{(ji)}e_{(ik)}=0\] by
relation (iii). Next, we note that $\sum_i \theta_i=2\theta$.
Hence, $\extd\alpha_i=\theta\alpha_i+\alpha_i\theta
=-\h\sum_k\alpha_k\alpha_i+\alpha_i\alpha_k=-\alpha_i^2$ as
required. \eproof

In the algebra $\CE_N$ the similar elements
\eqn{Etheta}{\theta_i=\sum_{i<j}[ij]-\sum_{j<i}[ji]=\sum_{j\ne i}
[ij]} form a commutative subalgebra isomorphic to the cohomology
of the flag variety\cite{FomKir:qua}. In our case we see that they
anticommute rather than commute. Also, while the elementary
symmetric polynomials of the $\theta_i$ in $\CE_N$ vanish, we have
in $\Lambda_N$
\[ \sum_i \theta_i=2\theta,\quad \sum_i\theta_i^2=0\]
as above. On the other hand, in \cite{FomKir:qua} the generators
$\theta_i$ are motivated from Dunkl operators on the cohomology of
the flag variety but in our case they have a direct noncommutative
geometrical interpretation as flat connections. We will see in the
next section that they are precisely the flat connections with
constant coefficients of minimal support.

\section{General differentials and flat connections up to $S_5$}

So far we have studied only one natural conjugacy class. However,
our approach associates a similar exterior algebra for any
nontrivial conjugacy class in a finite group $G$. Moreover, since
our constructions are $G$-invariant, we will obtain
`geometrically' plenty of $G$-modules naturally associated to the
conjugacy class. The cohomology does not tend to be a very
interesting representation but the moduli of flat connections
turns out to be more nontrivial and we will see that for $S_N$ it
does yield interesting irreducible modules. We begin with some
remarks for general finite groups $G$ equipped with a choice of
nontrivial conjugacy class.

First of all, the space of connections is an affine space. We take
as `reference' the form $-\theta$. Then one may easily see that
the differences $\phi\equiv \alpha+\theta$ transform covariantly
as \eqn{phi}{\phi=\sum_a\phi^a e_a\mapsto u\phi
u^{-1}=\sum_a{u\over R_a(u)}\phi^ae_a.} Moreover, the curvature of
$\alpha$ is \eqn{F}{
F(\alpha)=\extd\alpha+\alpha^2=\extd(\phi-\theta)
+(\phi-\theta)^2=\phi^2} in view of the properties of $\theta$.

\begin{lemma} Let $\alpha\in \Omega^1(G)$ be a connection. We
define its `cardinality' to be the number of nonzero components of
$\alpha+\theta$ in the basis $\{e_a\}$. This is gauge-invariant
and stratifies the moduli of connections.
\end{lemma} \proof A gauge transformation $u$ is invertible
hence the support of each component $\phi^a$ is gauge-invariant
under the transformation shown above. In particular, the number of
$\phi^a$ with nontrivial support is invariant. \eproof

We are particularly interested in invariant forms $\Lambda_{quad}$
and hence connections with constant coefficients $\phi^a$
(otherwise we have more refined gauge-invariant support data,
namely an integer-valued vector whose entries are the cardinality
of the support of each $\phi^a$). To simplify the problem further
we restrict to constant coefficients in $\{0,1\}$. We can project
any connection  to such a $\{0,1\}$-connection by replacing
non-zero $\phi^a$ by 1, so this limited class of connections gives
useful information about any connection.

\begin{propos} Flat connections with constant coefficients in
$\{0,1\}$ are in correspondence with subsets $X\subseteq \CC$ such
that
\[ \ad_x(X)=X,\quad \forall x\in X.\]
The intersection of such subsets provides a product in the moduli
of such flat connections that non-strictly lowers cardinality. The
stratum $F_n$ of subsets of a given cardinality $n$ is
$G$-invariant under $\Ad$.
\end{propos}
\proof The correspondence between $\{0,1\}$-connections and
subsets is via the support of the components $\phi^a$ regarded as
a function of $a\in \CC$ (so the cardinality of the connection is
that of the subset.) We have to solve the equation $\phi^2=0$. But
the relations in $\Lambda_{quad}$ are defined by the braiding
$\Psi$ and hence this equation is
\[ \Psi(\phi\tens\phi)=\phi\tens\phi.\]
Using the form of $\Psi$ this is
\[0= \sum_{a,b\in\CC} \phi^a\phi^b( e_{aba^{-1}}\tens e_a-e_a\tens
e_b)\] or
\[ \phi^a(\phi^{a^{-1}ba}-\phi^b)=0\quad\forall a,b\in \CC.\]
This translates into the characterisation shown. On the other
hand, this characterisation is clearly closed under intersection.
Finally, if $X$ is such a subset then $Y=\ad_g(X)$ is another such
subset because if $y=gxg^{-1}$ and $z=gwg^{-1}$ for $x,w\in X$
then $\ad_y(z)=(gxg^{-1})gwg^{-1}(gx^{-1}g^{-1})=g(\Ad_x(w)g^{-1}$
is in $Y$. \eproof

Over $\{0,1\}$ the stratum of top cardinality has one point,
$X=\CC$, which corresponds to $\alpha=0$ or $\phi=\theta$. The
stratum of zero cardinality likewise has one point, $X=\emptyset$
corresponding to $\alpha=-\theta$ or $\phi=0$, and the stratum
with cardinality $1$ can be identified with $\CC$, with
$\alpha=e_a-\theta$ or $\phi=e_a$ for $a\in \CC$. In between
these, the spans $kF_n$ are natural sources of permutation
$G$-modules. More precisely one typically has $\theta\in kF_n$ and
we look at the module \eqn{Vn}{ V_n=k F_n/k\theta.} For example,
the `submaximal stratum' (the one below the top one) associates a
representation to a conjugacy class, i.e. is an example of an
`orbit method' for finite groups. Like the usual orbit method for
Lie groups, it does not always yield an irreducible
representation, but does sometimes. It should be stressed that
this is only one example of the use of our geometrical methods to
define representations and we present it only as a first idea
towards a more convincing orbit method.

We now use the permutation groups $S_N$ to explore these ideas
concretely. For $S_3$ the full moduli of (unitary) flat
connections has already been found in \cite{MaRai:ele} for the
2-cycles class, while the other class is more trivial. The
exterior algebra in the second case is
$\Lambda_w=\Lambda_{quad}=k\<e_{123},e_{132}\>$ modulo the
relations
\[ e_{123}^2=0,\quad e_{132}^2=0,\quad e_{123}e_{132}
+e_{132}e_{123}=0\] (a Grassmann 2-plane). For brevity, we
suppress the brackets, so $e_{123}\equiv e_{(123)}$, etc.

\begin{table}
\[\begin{array}{l|c|l|c}
S_3& |\CC| & {\rm Solutions}/k\qquad \qquad{\rm Repn}\
kF_n/k\theta &{\rm Specht}\\ \hline (12)& 3 &
\begin{array}{lc}F_3=\{\theta\}& \\
F_1=\{e_{23},e_{13},e_{12}\}& \qquad{\rm fund}\\
\end{array}& {\rm fund}\\ \hline
(123)& 2& \begin{array}{lc}F_2=\{\cdot e_{123}+\cdot e_{132}\}& \\
F_1=\{e_{123},e_{321}\}&{}\quad\qquad{\rm sign}\\ \end{array}&
  {\rm sign}\\
\end{array}\]
\caption{Flat connections with constant coefficients on $S_3$ for
each conjugacy class, listed by cardinality. $\cdot$ denotes
independent nonzero multiples are allowed.} \end{table}

\begin{propos} The set of flat connections with
constant coefficients for $S_3,S_4$ with their various conjugacy
classes are as shown in Tables~2,3. For each stratum $F_n$ of
cardinality $n$, we list the corresponding $\phi$ up to an overall
scale. The associated representations $V_n$ turn out to be
irreducible.
\end{propos}
\proof This is done by direct computation. Note that the entries
$\phi$ of a discrete stratum each define a line of flat
connections $\alpha=\lambda\phi-\theta$ for a parameter $\lambda$.
The entries $\cdot e_{123}+\cdot e_{132}$, etc., define a plane of
connections $\alpha=\lambda e_{123}+\mu e_{132}$. As above, we
omit the brackets on the cycles labelling the $e_a$, for example
$e_{12,34}$ denotes $e_{(12)(34)}$. For the $V_n$ we enumerate the
flat connections in the stratum with coefficients $\{0,1\}$. The
resulting representation is then recognised using character
theory. Representations are labelled by dimension and by $\bar{\
}$ if the character at $(12)$ is negative. The fundamental
representation of $S_4$ means the standard 3-dimensional one.
\eproof

\begin{table}
\[\begin{array}{l|c|l|c}
S_4& |\CC| & {\rm
Solutions}/k\qquad\qquad\qquad\qquad\qquad\qquad\qquad\qquad\qquad
{\rm Repn}\ k F_n/k\theta&{\rm Specht}\\
 \hline (12)& 6 &
\begin{array}{lc}F_6=\{\theta\}& \\
F_3=\{\theta-\theta_i\}&{}\qquad\qquad\qquad\qquad{\rm fund}\\
F_2=\{\cdot e_{14}+\cdot e_{23},\cdot e_{13}+\cdot e_{24}, \cdot
e_{12}+\cdot e_{34}\}&{}\qquad\qquad\qquad\qquad 2\\
F_1=\{e_{a}\}& \\
\end{array}&{\rm fund}\\
\hline (12)(34)& 3& \begin{array}{lc}F_3=\{\cdot e_{12,34}+\cdot
e_{13,24}+\cdot e_{14,23}\}& \\
F_2=\{\cdot e_{12,34}+\cdot e_{13,24},\cdot e_{13,24}+\cdot
e_{14,23},\cdot e_{12,34}+\cdot e_{14,23}\}&\qquad\quad 2\\
F_1=\{e_{a}\}& \end{array}& 2\\
\hline (123)& 8 & \begin{array}{lc}F_8=\{\cdot( e_{123}+e_{142}
+e_{134}+e_{243})
+\cdot(e_{132}+e_{124}+e_{143}+e_{234} )\}& \\
F_4=\left\{\begin{array}{c}e_{123}+e_{142}+e_{134}+e_{243},\\
e_{132}+e_{124}+e_{143}+e_{234}\\ \end{array}\right\}
&\kern -.5in {\rm sign}\\
F_2=\{\cdot e_{123}+\cdot e_{132},\cdot e_{142}+\cdot
e_{124},\cdot e_{134}+\cdot e_{143},\cdot
e_{243}+\cdot e_{234}\}&\kern -.5in {\rm fund} \\
F_1=\{e_{a}\}&\\
\end{array}& \overline{\rm fund}\\
\hline (1234)& 6& \begin{array}{lc} F_6=\{\theta\}& \\
F_2=\{\cdot e_{1234}+\cdot e_{1432}, \cdot e_{1243}+\cdot
e_{1342},\cdot e_{1324}+\cdot e_{1423}\}&\qquad\qquad2\\
F_1=\{e_{a}\}& \\ \end{array}& {\rm sign}\\
\end{array}\]
\caption{Flat connections with constant coefficients on $S_4$ for
each conjugacy class, listed by cardinality. $\cdot$ denotes
independent nonzero multiples are allowed.} \end{table}

For comparison, the tables also list the standard Specht module of
the Young tableau of conjugate shape to that of the conjugacy
class. We see that our `orbit method' produces comparable
(although different) answers. As was to be expected, we do not
obtain all irreducibles from consideration of $\{0,1\}$
connections alone. Similarly for the $S_5$ case:

\begin{table}
\[\begin{array}{l|c|l|c}
S_5& |\CC| & {\rm Solutions}/\Z_2 &
{\rm Repn}\ k F_n/k\theta\\
 \hline (12)& 10 &
\begin{array}{l}
F_{10}=\{\theta\}\\ F_6=\{\theta-\theta_i\}\\
| F_4| =10,\ |F_3|=10,\ |F_2|=15\\
F_1=\{e_a\}\\
\end{array}&{\rm fund}\\
\hline (12)(34)& 15 &
\begin{array}{l}
F_{15}=\{\theta\}\\
F_5=\left\{\begin{array}{c}
e_{14,23}+e_{12,35}+e_{13,45}+e_{25,34}+e_{15,24},\\
e_{14,23}+e_{13,25}+e_{24,35}+e_{15,34}+e_{12,45},\\
e_{13,24}+e_{12,35}+e_{24,35}+e_{12,45}+e_{15,24},\\
e_{13,24}+e_{15,23}+e_{25,34}+e_{12,45}+e_{14,35},\\
e_{12,34}+e_{13,25}+e_{23,45}+e_{15,24}+e_{14,35},\\
e_{12,34}+e_{15,23}+e_{24,35}+e_{14,25}+e_{13,45}
\\ \end{array}\right\} \\
|F_3|=15, \ |F_2|=15\\
F_1=\{e_a\}\\ \end{array}& \bar 5\\
\hline (123)& 20&\begin{array}{l} F_{20}=\{\theta\}\\
F_8=\{\theta-\theta_i\}\\
|F_4|=10,\ |F_2|=10\\
F_1=\{e_a\}\\ \end{array}& {\rm fund}\\
\hline (123)(45)& 20& \begin{array}{l}
F_{20}=\{\theta\}\\
F_2=\{e_{xyz,12}+e_{xzy,12},\cdots ({\rm
all\ 2-cycles;\ xyz\ {\rm complementary}})\}\\
F_1=\{e_a\}\\
\end{array}& {\rm fund}\oplus 5\\
\hline (1234)& 30 &
\begin{array}{l} F_{30}=\{\theta\}\\
F_{10}=\left\{{\scriptscriptstyle\begin{array}{c}
e_{1234}+e_{1523}+e_{2435}+e_{2534}+e_{1245}\\
{\ }\quad +e_{1542}+e_{1354}+e_{1453}+e_{1432}+e_{1325},\\
e_{1234}+e_{1253}+e_{2453}+e_{2354}+e_{1524}\\
{\ }\quad +e_{1425}+e_{1345}+e_{1543}+e_{1432}+e_{1352},\\
e_{1234}+e_{1253}+e_{2453}+e_{2354}+e_{1524}\\
{\ }\quad +e_{1425}+e_{1345}+e_{1543}+e_{1432}+e_{1352},\\
e_{1523}+e_{2453}+e_{2354}+e_{1243}+e_{1452}\\
{\ }\quad +e_{1254}+e_{1534}+e_{1435}+e_{1342}+e_{1325},\\
e_{1532}+e_{2435}+e_{2534}+e_{1452}+e_{1254}\\
{\ }\quad +e_{1345}+e_{1324}+e_{1543}+e_{1423}+e_{1235},\\
e_{1253}+e_{2345}+e_{2543}+e_{1245}+e_{1542}\\
{\ }\quad +e_{1534}+e_{1435}+e_{1324}+e_{1423}+e_{1352}\\
\end{array}}\right\}\\
F_6=\{\theta-\theta_i\}\\
|F_5|=12,\ |F_2|=15\\
F_1=\{e_a\}\\ \end{array}& \bar 5\\
\hline (12345)& 24 &\begin{array}{l} F_{24}=\{\theta\}\\
F_{12}=\left\{\begin{array}{c}
e_{12345}+e_{12453}+e_{12534}+\cdots ({\rm
sum\ even}),\\
  e_{12354}+e_{12435}+e_{12543}+\cdots ({\rm sum\ odd})\\
  \end{array}\right\}\\
|F_4|=6,\ |F_3|=24,\ |F_2|=36 \\
F_1=\{e_a\}\\ \end{array}& {\rm sign}\\
\end{array}\]
\caption{Flat connections with constant coefficients in $\{0,1\}$
on $S_5$ for each conjugacy class. The submaximal strata are
listed in detail, as well as the associated representation.}
\end{table}

\begin{propos} The set of
flat connections with constant coefficients $\{0,1\}$ for $S_5$
with its various conjugacy classes are as shown in  Table~4,
organised by stratum $F_n$ of cardinality $n$. The submaximal
strata are shown in detail as well as their associated
representations $V_n$.
\end{propos}
\proof These results have been obtained with GAP to compute $\ad$
tables followed by MATHEMATICA running for several days to
enumerate the flat connections. In the tables $\theta_i$ denotes a
sum over the relevant size cycles containing $i$, extending our
previous notation. \eproof

We see that for low $N$ the $V_n$ tend to be irreducible, but they
are not always. Notably, the $(123)(45)$ conjugacy class for $S_5$
has a 9-dimensional representation associated to the submaximal
stratum. A possible refinement would be to consider only the
`discrete series' i.e. flat connections where $\phi$ of a fixed
normalisaiton is not deformable. For $S_4$ this means from Table~3
the strata $F_3$ for the 2-cycles class and $F_4$ for 3-cycles. A
different problem is that one does not get all irreducibles in
this way (since one only gets permutation modules). To go beyond
this one could consider the full moduli of flat connections
including the non-discrete series but with other constraints. For
example, one could consider connections with values in
$\{-1,0,1\}$, or one could introduce further geometric ideas such
as `polarizations' to our noncommutative setting.

Finally, each of our conjugacy classes on $S_N$ has its associated
quadratic algebra $\Lambda_{quad}$ of interest in its own right.
We give just one example.

\begin{propos} For $S_4$ with its 3-cycle conjugacy class
the exterior algebra $\Lambda_{quad}$ has relations \begin{align*}
&e_{xyz}^2=0,\quad e_{xyz}
e_{xzy}+e_{xzy}e_{xyz}=0,\\
 &e_{123}e_{134}+e_{134}e_{142}+e_{142}e_{123}=0,\quad
e_{123}e_{243}+e_{243}e_{134}+e_{134}e_{123}=0,\\
& e_{134}e_{243}+e_{243}e_{142}+e_{142}e_{134}=0,\quad
e_{123}e_{142}+e_{142}e_{243}+e_{243}e_{123}=0,\\
 &e_{123}e_{124}+e_{124}e_{134}+e_{134}e_{234}+e_{234}e_{123}=0,\\
&e_{123}e_{143}+e_{143}e_{243}+e_{243}e_{124}+e_{124}e_{123}=0,\\
 &e_{123}e_{234}+e_{234}e_{142}+e_{142}e_{143}+e_{143}e_{123}=0\end{align*}
and their seven conjugate-transposes (i.e. replacing $e_{xyz}$ by
$e_{xzy}$ and reversing products).
\end{propos}
\proof Direct computation of the kernel of $\id-\Psi$ using GAP
and MATHEMATICA. We omit the brackets around the 3-cycle labels
(as above). \eproof

We conclude with one general result pertaining to the above ideas.

\begin{propos} For the 2-cycles conjugacy class on $S_N$, the flat
connections with constant coefficients in $\{0,1\}$ of submaximal
cardinality are precisely the $\alpha_i=-\theta_i$ in
Proposition~3.4. The associated module is the fundamental
representation of $S_N$.
\end{propos}
\proof The $\alpha_i$ correspond to $\phi_i=\theta-\theta_i$ and
have cardinality $({N-1\atop 2})$. Consider any flat connection
with constant coefficients in $\{0,1\}$ with corresponding $\phi$
or corresponding subset $X$  in Proposition~4.2 of cardinality
$|X|\ge ({N-1\atop 2})$. Suppose there exists $i\in\{1,\cdots,N\}$
such that for all $i'$, $(ii')\notin X$. But there are only
$({N-1\atop 2})$ such elements of $\CC$ (those not containing $i$
in the 2-cycle) so $X$ cannot have cardinality greater than this,
hence $|X|=({N-1\atop 2})$ and $\phi=\phi_i$. Otherwise, we
suppose that for all $i$ there exists $i'$ such that $(ii')\in X$.
Then for any $i,j$ we have $(ii'),(jj')\in X$ hence by the $\Ad$
closure of $X$ we have $(ij)\in X$, i.e. $X=\CC$ or $\phi=\theta$.
\eproof

\section{Metric structure}

In this section we look at some more advanced aspects of the
differential geometry for $S_N$, but for the 2-cycle calculus.
First of all, just as the dual of the invariant 1-forms on a Lie
group can be identified with the Lie algebra, the space $\CL=
\Lambda^{1*}$ for a bicovariant differential calculus on a
coquasitriangular  Hopf algebra $A$ is typically a braided-Lie
algebra in the sense introduced in \cite{Ma:lie}. Moreover, every
braided-Lie algebra has an enveloping algebra\cite{Ma:lie} which
in our case means \eqn{UL}{
U(\CL)=T\Lambda^{1*}/\image(\id-\Psi^*)= \Lambda_{quad}^!,} where
$!$ is the quadratic algebra duality operation. There is also a
canonical algebra homomorphism $U(\CL)\to H$ where $H$ is dual to
$A$. We will call a differential structure `connected' if this is
a surjection. This theory applies to the Drinfeld-Jimbo $U_q(g)$
and gives it as generated by a braided-Lie algebra for each
connected calculus.

However, the theory also applies to finite groups and in this case
the axioms of a braided-Lie algebra reduce to what is called in
algebraic topology a rack. Thus, given a conjugacy class on a
finite group $G$, the associated rack or braided-Lie algebra
is\cite{Ma:rief} \eqn{LG}{ \CL=\{x_a\}_{a\in \CC},\quad
[x_a,x_b]=x_{b^{-1}ab},\quad \Delta x_a=x_a\tens x_a,\quad
\eps(x_a)=1.} The analogue of the Jacobi identity is \eqn{JacG}{
[[x_a,x_c],[x_b,x_c]]=[[x_a,x_b],x_c].} The enveloping algebra is
the ordinary bialgebra $U(\CL)=k\<x_a\>$ modulo the relations
$x_ax_b=x_b x_{b^{-1}ab}$ and its homomorphism to the group
algebra of $G$ is $x_a\mapsto a$. This is surjective precisely
when any element of $G$ can be expressed as a product of elements
of $\CC$, i.e. by a path with respect to our differential
structure (which determines the allowed steps as elements of
$\CC$) connecting the element to the group identity. Thus, in our
finite group setting, the quadratic algebra $\Lambda_{quad}$ is
the $!$-dual of a fairly natural quadratic extension of the group
algebra as an infinite-dimensional bialgebra. Note also that the
flat connections in Proposition~4.2 define braided sub-Lie
algebras.

Next, associated to any braided-Lie algebra is an $\ad$-invariant
and braided-symmetric (with respect to $\Psi$) braided-Killing
form, which may or may not be nondegenerate. This is computed in
\cite{Ma:rief} for finite groups and one has
\eqn{KillG}{\eta^{a,b}\equiv\eta(x_a,x_b)=\#\{c\in\CC|\
cab=abc\}.} The associated metric tensor in
$\Omega^1\tens_{k(G)}\Omega^1$ is
\[ \eta=\sum_{a,b}\eta^{a,b}e_a\tens e_b\]
It is easy to see that among $\ad$-invariant $\eta$, `braided
symmetric' under $\Psi$ is equivalent to symmetric in the usual
sense. It is also equivalent (by definition of $\wedge)$ to
$\wedge(\eta)=0$ under the exterior product.

\begin{propos} For the braided-Lie algebra associated to the
2-cycle calculus on $S_N$, the braided-Killing form is
\[ \eta^{(ij),(ij)}=\left({N\atop 2}\right),\quad
\eta^{(ij),(km)}=\left({N-4\atop 2}\right)+2,\quad
\eta^{(ij),(jk)}=\left({N-3\atop 2}\right)\] for distinct
$i,j,k,m$. Moreover, the calculus is `connected'.
\end{propos}
\proof All of $\CC$ commutes with $(ij)^2=e$. In the second case
all elements disjoint from $i,j,k,m$ and $(ij),(km)$ themselves
commute with $(ij)(km)$. For the third case all elements disjoint
from $i,j,k$ commute with $(ij)(jk)$.  The connectedness is the
well-known property that the 2-cycles can be taken as generators
of $S_N$. \eproof

To be a metric, we need $\eta$ to be invertible, which we have
verified explicitly at least up to $N<30$. Other symmetric and
invariant metrics also exist, not least
$\eta^{a,b}=\delta_{a,b^{-1}}$, the Kronecker $\delta$-function
which is always invertible and works for any conjugacy class on
any finite group that is stable under inversion. The general
situation for $S_N$ is:

\begin{propos} The most general conjugation-invariant metric for the 2-cycle
calculus on $S_N$ has the symmetric form
\[  \eta^{(ij),(ij)}=\alpha,\quad \eta^{(ij),(km)}=\beta,\quad
\eta^{(ij),(jk)}=\gamma\] for distinct $i,j,k,m$, where
$\alpha,\beta,\gamma$ are three arbitrary constants. Moreover,
\[ \det(\eta)=\left(\alpha+\beta-2\gamma\right)^{N(N-3)\over 2}
\left(\alpha-(N-3)\beta+(N-4)\gamma\right)^{N-1}\]
\[ {\ }\quad\quad\quad\cdot\left(\alpha+{(N-2)(N-3)\over 2}
\beta+2(N-2)\gamma\right)\] at
least up to $N\le 10$.
\end{propos}
\proof Invariance here means $\eta^{gag^{-1},gbg^{-1}}=\eta^{a,b}$
for all $g\in G$. We use the mutually exclusive  notations $a=b$,
$a\perp b$ and $a\sim b$ as in the proof of Theorem~3.2, which is
clearly an $\ad$-invariant decomposition of $\CC\times\CC$  (since
the action of $S_N$ is by a permutation of the 2-cycle entries).
Clearly all the diagonal cases $a=b$ have the same value since
$\CC$ is a conjugacy class. Moreover, any $(ij)\perp (km)$ (for
$N\ge 4$) is conjugate to $(12)\perp (34)$ by the choice of a
suitable permutation (which we use to make the conjugation), so
all of these have the same value. Similarly every $(ij)\sim(jk)$
(for $N\ge 3$) is conjugate to $(12)\sim(23)$, so these all have
the same value. We then compute the determinants for $N\le 10$ and
find that they factorise in the form stated. The first two factors
cancel in the case of $N=2$.  \eproof

Armed with an invertible metric, one may compute the associated
Hodge-* operator, etc. as in \cite{MaRai:ele} for $S_3$. The
computation of this for $S_N$ is beyond our present scope as it
would require knowledge of $\Lambda_{quad}$ in all degrees (we do
not even know the dimensions for large $N$). It is also beyond our
scope to recall all the details of noncommutative Riemannian
geometry, but along the same lines as for $S_3$ in \cite{Ma:rief}
we would expect a natural regular Levi-Civita connection with
Ricci curvature tensor proportional to the metric modulo
$\theta\tens\theta$. Moreover, the same questions can be examined
for the other conjugacy classes or `Riemannian manifold'
structures on $S_N$.

\section{Braided group structure on $\CE_N$}

In this section we show that that the Fomin-Kirillov algebra
$\CE_N$ is a Hopf algebra in the braided category of crossed
$S_N$-modules. In fact, we will find that like the exterior
algebras $\Lambda_N$, it is a `braided linear space' with additive
coproduct on the generators\cite{Ma:book}. We recall that a
braided group $B$ has a coproduct $\und\Delta:B\to B\und\tens B$
which is coassociative and an algebra homomorphism provided the
algebra $B\und\tens B$ is the braided-tensor product where
\[ (a\tens b)(c\tens d)=a\Psi_{B,B}(b\tens c)d\] where $a,b,c,d\in
B$ and $\Psi_{B,B}$ is the braiding on $B$. We show how the cross
product (usual) Hopf algebras $kS_N\rbiprod \CE_N$ in
\cite{FomPro:fib} and the skew derivations $\Delta_{ij}$ related
to divided differences in \cite{FomKir:qua} arise immediately as
corollaries of the braided group structure. While the $\CE_N$ are
already-well studied by explicit means, we provide a more
conceptual approach that is also more general and applies both to
other conjugacy classes and to other groups beyond $S_N$.

As in \cite{FomKir:qua} we consider that the algebra $\CE_N$ is
generated by an $[{N\atop 2}]$-dimensional vector space $E_N$
(say) with basis $[ij]$ where $i<j$, and we extend the notation to
$i>j$ by $[ij]=-[ji]$.

\begin{theorem} The algebras $\CE_N$ are `braided groups' or
Hopf algebras in the category of $S_N$-crossed modules. Here
\[ g.[ij]=[g(i)\ g(j)]=\begin{cases}[g(i)\ g(j)]& {\rm if}\
g(i)<g(j)\\ -[g(j)\ g(i)]&{\rm if}\ g(i)>g(j)\end{cases},\ \forall
g\in S_N,\quad |[ij]|=(ij)\] is the crossed module structure on
$E_N$, where $|\ |$ denotes the $S_N$-degree. Let $\Psi$ denote
the induced braiding, then
\[ \CE_N=T E_N/\ker(\id+\Psi),\quad\und\Delta [ij]=[ij]\tens 1+1\tens [ij],
\quad \und\eps[ij]=0\] is an additive braided group or `linear
braided space' in the category of $S_N$-crossed modules.
\end{theorem}
\proof It is easy to verify that this is a crossed module
structure. Thus $|g.[ij]|=\pm (g(i)\ g(j))=\pm g(ij)g^{-1}=
|g[ij]g^{-1}|$ for the two cases (note that we consider the
$S_N$-degree extended by linearity). The braiding is then \[
\Psi([ij]\tens [km])=(ij).[km]\tens [ij]\] as defined by the
crossed module structure. This is a signed version of the braiding
used in Proposition~3.1 and by a similar analysis to the proof
there, one finds that the kernel of $\id+\Psi$ is precisely
spanned by the relations of $\CE_N$. In particular, note that
\[ \Psi([ij]\tens [ij])=-[ij]\tens[ij],\quad \Psi([ij]\tens [km])
=[km]\tens[ij]\] if disjoint, which gives the relations $[ij][ij]
=0$ and $[ij][km]=[km][ij]$ when disjoint. Similarly for the
3-term relations $[ij][jk]+[jk][ki]+[ki][ij]=0$ when $i,j,k$ are
distinct. Next, we define the coalgebra structure on the
generators as stated and verify that these extend in a
well-defined manner to a braided group structure on $\CE_N$.
there. This part is the same as for any braided-linear space
\cite{Ma:book} and we do not repeat it. The only presentational
difference is that we directly define the relations as
$\ker(\id+\Psi)=0$ rather then seeking some other matrix $\Psi'$
such that $\image(\id-\Psi')=\ker(\id+\Psi)$. \eproof

\begin{corol} If $\CE_N$ is finite-dimensional then it has a unique element
in top degree. \end{corol} \proof A top degree element would be an
integral in the braided-Hopf algebra. But as for a usual Hopf
algebra, the integral if it exists is unique (a formal proof in
the braided case is in \cite{Lyu:mod}\cite{BKLT:int}). \eproof

Also, the biproduct bosonisation of any braided group $B$ in the
category of left $A$-crossed modules is an ordinary Hopf algebra
$B\lbiprod A$ (where $A$ is an ordinary Hopf algebra with
bijective antipode). This is the simultaneous cross product and
cross coproduct in the construction of \cite{Rad:str}, in the
braided group formulation \cite[Appendix]{Ma:skl}. In our case $A$
is finite dimensional so $B$ also lives in the category of right
$A^*$-crossed modules. Hence we immediately have two ordinary Hopf
algebras, the first of which recovers the cross product observed
in \cite{FomPro:fib} and studied further there.

\begin{corol} Biproduct bosonisation in the category of left
crossed $S_N$-module structure gives an ordinary Hopf algebra
$\CE_N\lbiprod kS_N$ with \[ g[ij]=[g(i)\, g(j)]g,\quad\forall
g\in S_N,\quad \Delta [ij]=[ij]\tens 1+(ij)\tens[ij],\quad
\eps[ij]=0\] extending that of $kS_N$, as in \cite{FomPro:fib}.
Bosonisation in the equivalent category of right $k(S_N)$-crossed
modules gives an ordinary Hopf algebra $k(S_N)\rbiprod \CE_N$ with
\[ [ij]f=R_{(ij)}(f)[ij],\quad\forall f\in k(S_N),\quad
\Delta[ij]=\sum_{g\in S_N}[g(i)\, g(j)]\tens \delta_g+1\tens
[ij].\]
\end{corol}
\proof The $kS_N$-module structure defines the cross product and
the $kS_N$-coaction $\Delta_L[ij]=(ij)\tens [ij]$ defines the
cross coproduct. In the second case the $S_N$-grading defines an
action of $k(S_N)$ and the $kS_N$-module structure defines the
$k(S_N)$-coaction $\Delta_R[ij]=\sum_g[g(i)\, g(j)]\tens (ij)$ by
dualisation. \eproof

Next, from a geometrical point of view the $\CE_N$ are `linear
braided spaces', i.e. the coproduct $\und\Delta$  corresponds to
the additive group law on usual affine space in terms of its usual
commutative polynomial algebra in several variables, but now in a
braided-commutative version. We will use several results from this
theory of linear braided spaces. For clarity we explicitly label
the generators of $\CE_N$ by 2-cycles. Thus $[ij]=e_{(ij)}$ if
$i<j$. The products are different from those of $\Lambda_N$ but we
identify the basis of generators. In this notation we have
\eqn{zetaG}{ g.e_b=\zeta_{g,b}e_{gbg^{-1}},\quad |e_a|=a}
\eqn{zeta}{ \zeta_{(ij),(ij)}=-1,\quad \zeta_{(ij),(km)}=1,\quad
\zeta_{(ij),(jk)} =\begin{cases}
1&{\rm if}\ i<j<k\\ 1&{\rm if}\ j<i<k\\
-1&{\rm if}\ j<k<i\\ -1&{\rm if}\ i<k<j\\ 1&{\rm if}\ k<i<j\\
1&{\rm if}\ k<j<i\end{cases}} for $i,j,k,m$ distinct, where
$\zeta$ extends to $S_N$ in its first argument by
$\zeta_{gh,b}=\zeta_{g,hbh^{-1}}\zeta_{h,b}$ for all $g,h\in S_N$,
i.e.,
\[\zeta\in Z^1_{\ad}(S_N,k(\CC));\quad \zeta(g)(b)=\zeta_{g,b},\]
as a multiplicative cocycle (using the multiplication of the
algebra $k(\CC)$ of functions on $\CC$ and with $\ad$ the right
action on $k(\CC)$ induced by conjugation). Thus, the algebras
$\CE_N$ differ from the exterior algebras $\Lambda_N$ precisely by
the introduction of a cocycle. This makes precise how to construct
analogues of the $\CE_N$ for other finite groups.

In this notation we have for any braided linear
space\cite{Ma:fre}\cite{Ma:book} \eqn{Deltaprod}{
\und\Delta(e_{b_1}\cdots e_{b_m})=\sum_{r=1}^m e_{c_1}\cdots
e_{c_r}\tens e_{c_{r+1}}\cdots e_{c_m}\left[{m\atop
r},\Psi\right]^{c_1\cdots c_m}_{b_1\cdots b_m}} on products of
generators. We view the braiding $\Psi$ as a matrix  (denoted $PR$
in\cite{Ma:book}). The braided binomial matrices have been
introduced by the author in exactly this context and are not
assumed to be invertible. There is also a braided antipode $\und
S:\CE_N\to\CE_N$ defined as $-1$ on the generators and extended
braided-antimultiplicatively in the sense $\und
S(fg)=\cdot\Psi(\und S f\tens\und S g)$ for all $f,g\in\CE_N$, see
\cite{Ma:book}. In our case this comes out inductively as
\eqn{antmult}{ \und S(e_a f)=-\cdot\Psi(e_a\tens \und S f)
=-(a.\und S f )e_a,\quad \forall f\in \CE_N.}

Next we note that the braided group $\CE_N$ is certainly
finite-dimensional in each degree, so it has a graded-dual braided
group $\CE_N^*$. However, a linear braided space and its dual can
typically be identified in the presence of an invariant metric. In
our case we use the Kronecker $\delta_{a,b}$ metric and have:

\begin{propos} $\CE_N$ is self-dually paired as a braided group, with
pairing
\[\<e_{a_n}\cdots e_{a_1},e_{b_1}\cdots e_{b_m}\>= \delta_{n,m}
([n,\Psi]!)^{a_1\cdots a_n}_{b_1\cdots b_n}\]
\end{propos}
\proof The pairing we take on the generators is
$\<e_a,e_b\>=\delta_{a,b}$, which is compatible with the
$S_N$-grading since all $a$ have order 2, and compatible with the
action of $g\in S_N$ since $(\pm 1)^2=1$, i.e. the pairing is a
morphism to the trivial crossed module. We extend this to products
via the axioms of a braided group as explained in \cite{Ma:book}
to obtain the pairing stated, using the above formula for
$\und\Delta$ and properties of the braided binomial operators in
relation to braided factorial matrices. It follows from the
construction that the pairing is well-defined in its second input.
This part is the same as in \cite{Ma:book}. It is also
well-defined in its first input after we observe that
$\Psi^*=\Psi$, where $\Psi^*$ is defined as the adjoint on
$E_N\tens E_N$ with respect to the braided-tensor pairing
$E_N\tens E_N\tens E_N\tens E_N\to k$ (in which we apply $\<\ ,\
\>$ to the inner $E_N\tens E_N$ first and then the outer two.)
\eproof

This implies in particular that the two Hopf algebras
$k(S_N)\rbiprod \CE_N$ and $\CE_N\lbiprod kS_N$ in Corollary~6.3
are dually paired. There are many more applications of the
braided-linear space structure. As a less obvious one we compute
the braided-Fourier theory \cite{KemMa:alg} introduced for
q-analysis on braided spaces.

\begin{propos} For $\CE_3$ the coevaluation for the pairing in
Proposition~6.4 is \begin{align*} \exp &=1\tens
1+[12]\tens[12]+[23]\tens[23]+[31]\tens[31]\\
&-[12][23]\tens[12][31]+[23][12]\tens[12][23]+[23][31]\tens[31][23]
-[31][23]\tens[31][12]\\
&+[31][12][23]\tens[31][12][23]+[12][23][31]\tens[12][23][31]
+[23][31][12]\tens[23][31][12]\\
&+[12][23][12][31]\tens[12][23][12][31]\end{align*}
 and this along
with the integration $\int$ defined as the coefficient of the top
element $[12][23][12][31]$ (and zero in lower degree) defines
braided Fourier transform $\CS:\CE_3\to \CE_3$
\begin{align*} &\CS(1)=[12][23][12][31],\quad
\CS([12])=[31][12][23],\quad
\CS([23])=[12][23][31]\\
& \CS([31])=[23][31][12],\quad \CS([12][23])=[31][12],\quad
\CS([23][12])=[31][23]\\
&\CS([23][31])=[12][23],\quad \CS([31][23])=[12][31],\quad
\CS([31][12][23])=-[12]\\
&\CS([12][23][31])=-[23],\quad \CS([23][31][12])=-[31],\quad
\CS([12][23][12][31])=1.\end{align*} It obeys $\CS^2=\id$ in
degrees 0,4, $\CS^2=-\id$ in degrees 1,3 and $\CS^3=\id$ in degree
2.
\end{propos}
\proof Let $\{e^{(r)}_A\}$ be a basis of $\CE_N$ in degree $r$ and
$\{f^{(r)A}\}$ the dual basis with respect to the pairing. In the
nondegenerate case this is given by the inverse of the quotient
operator $[r,\Psi]!$ acting on the degree $r$ component of
$\CE_N$. The coevaluation for the pairing is
\[ \exp=\sum_r \sum_A e_A^{(r)}\tens f^{(r)A}\]
which computes as stated for $\CE_3$.  We take basis
$\{[12],[23],[31]\}$ for degree 1, which is orthonormal with
respect to the duality pairing. For degree 2 we take basis
$\{[12][23],[23][12],[23][31],[31][23]\}$ with dual
$\{-[12][31],[12][23],[31][23],-[31][12]\}$. The basis in degree 3
is the Fourier transform of the basis in degree 1 and orthonormal.
The braided Fourier transform is defined on general braided groups
possessing duals and integrals \cite{KemMa:alg}. Here
\[ \CS(f)=(\int\tens\id) f \exp,\quad \forall f\in \CE_N\]
which computes as stated using the relations of $\CE_3$. A similar
formula including the braided antipode $\und S$ provides the
inverse Fourier transform. We compute $\CS^2$ as stated. In fact
$\CS^2=\CT$, where $\CT(f)=|f|.f$ on $f\in \CE_N$ of homogeneous
$S_N$-degree. One also has $\CS=\CT^{-1}=-\und S$ in degree 2.
 \eproof

Also from the linear braided space structure, $\CE_N^*$ and in our
case $\CE_N$ acts on the algebra $\CE_N$ by infinitesimal
translation from the left and right, which means respectively
partial derivatives $D_a,\bar D_a:\CE_N\to \CE_N$ for each $a\in
\CC$ (they are denoted $\del^a,\bar\del^a$ in the general theory
of \cite{Ma:book}). Thus the left partial derivative is defined as
the coefficient of $e_a\tens$ in the operator $\und\Delta$, which
from the above yields \eqn{delprod}{ D_a(e_{a_1}\cdots
e_{a_m})=e_{b_2}\cdots e_{b_m}[m,\Psi]^{a b_2\cdots
b_m}_{a_1\cdots a_m}.} As for any braided linear space these
necessarily represent $\CE_N$ on itself and obey \eqn{Dleib}{
D_a(fg)=D_a(f)g+\cdot\Psi^{-1}(D_a\tens f)g,\quad \forall f,g\in
\CE_N,} making $\CE_N$ an opposite braided $\CE_N$-module algebra
(i.e., in the braided category with inverted braiding). Similarly
the right partial derivatives $\bar D_a$ are defined via right
translations but converted into an action from the left via the
braiding. This yields \eqn{bardelprod}{ \bar D_a(e_{a_1}\cdots
e_{a_m})=e_{b_2}\cdots e_{b_m}[m,\Psi^{-1}]^{a b_2\cdots
b_m}_{a_1\cdots a_m}} and necessarily represent $\CE_N$ on itself
with \eqn{bardelleib}{ \bar D_a(fg)=\bar D_a(f)g+\cdot \Psi(\bar
D_a\tens f)g,\quad \forall f,g\in\CE_N,} making $\CE_N$ a
$\CE_N$-module algebra in its original braided category. The
partial derivatives and their conjugates are related by the
braided antipode according to \eqn{Sdel}{ \und S D_a=-\bar D_a
\und S} as shown in \cite{Ma:qsta}. Proofs of all of these facts
are by braid-diagram methods as part of our established theory of
braided groups.

\begin{corol} In the case of $\CE_N$ the braided partial
derivatives $D_a,\bar D_a$ are covariant (morphisms in the
category of crossed modules) in the sense
\[ |D_a f|=a|f],\quad g.D_a (f)=\zeta_{g,a}
D_{gag^{-1}}(g.f),\quad \forall g\in S_N, \ f\in\CE_N\] (and
similarly for $\bar D_a$). They obey $D_a(e_b)=\bar
D_a(e_b)=\delta_{a,b}$ and the braided Leibniz rules
\[ D_a(fg)=D_a(f)g+f \zeta_{|f|^{-1},a} D_{|f|^{-1}a|f|}(g),
\quad \bar D_a(fg)=\bar D_a(f)g+(a.f)\bar D_a (g)\] for all
$f,g\in \CE_N$ and $f$ of homogeneous $S_N$-degree $|f|$ in the
first case.
\end{corol}
\proof Whereas the above review of $D_a,\bar D_a$ holds for any
additive braided group as part of a general theory, we specialize
now to the particular braided category for the case of $\CE_N$.
First of all, the $D_a,\bar D_a$ are defined above as evaluation
on $e_a$ of morphisms $D,\bar D:E_N\tens\CE_N\to \CE_N$. Thus
$D(g.e_a\tens g.f)=D(\zeta_{g,a}e_{gag^{-1}}\tens
g.f)=g.D(e_a\tens f)$ translates to the condition as shown for
$g\in S_N$ and $f\in \CE_N$. Likewise, commuting with the total
$S_N$-degree gives the other part of the morphism condition. Next
we compute the braiding and its inverse on $E_N\tens\CE_N$ in the
category of left crossed modules for the crossed module structure
stated. In general $\Psi(f\tens g)=|f|.g\tens f$ for $f$ of
homogeneous $S_N$-degree $|f|$. Applying $D,\bar D$ yields the
Leibniz rules as stated with
\[ \Psi^{-1}(D_a\tens f)=\zeta_{|f|^{-1},a}\, f\tens
D_{|f|^{-1}a|f|},\quad\Psi(\bar D_a\tens f) =a.f\tens \bar D_a.\]
The braiding of $D_a$ in these expressions corresponds by
definition to the braiding of the element $e_a$ which it
represents (similarly for $\bar D_a$).  \eproof

These $\bar D_{(ij)}$ therefore coincide with the operators
denoted $\Delta_{ij}$ in the notation of \cite{FomKir:qua}. It is
proven there that their restriction to polynomials in the
$\{\theta_i\}$ (i.e. to the cohomology of the flag variety) yields
the finite difference operators \eqn{divdel}{
\del_{ij}f={f-(ij).f\over \theta_i-\theta_j}} where $(ij).f$
interchanges the $i,j$ arguments of $f(\theta_1,\cdots,\theta_N)$.
We see that these $\bar D_a$ follow directly from the braided
group structure as infinitesimal translations, which ensures that
they are well-defined and form a representation of $\CE_N$ on
itself. It also provides computational tools, for example
braided-Fourier transform intertwines the braided derivatives with
multiplication in $\CE_N$, as shown in general in
\cite{KemMa:alg}. Another canonically-defined representation of a
braided group on itself with a similar braided-Leibniz property to
the $\bar D_a$ is the braided adjoint action, which comes out for
$\CE_N$ as \eqn{braad}{ \und\ad_{e_a}(f)\equiv e_a
f+\cdot\Psi(\und S e_a\tens f)=e_a f-(a.f)e_a} \eqn{adleib}{
\und\ad_{e_a}(fg)=\und\ad_{e_a}(f)g+(a.f)\und\ad_{e_a}(g).} This
has no direct geometrical analogue (the usual polynomial algebra
is commutative so that $\und\ad$ is zero).

Let us also note that the full bosonisation
theorem\cite[Thm.~9.4.12]{Ma:book} of $\CE_N$ provides another
ordinary Hopf algebra $\CE_N\lbiprod D(k S_N)$ such that its
category of modules is fully equivalent to the category of braided
modules of the braided group $\CE_N$. In particular,
module-algebras of this ordinary Hopf algebra are the same thing
as braided $\CE_N$-module algebras, such as provided by $\bar
D_a,\ad_{e_a}$ above. The Drinfeld double of a finite group is
itself a semidirect product $D(kS_N)=k(S_N){}_{\ad}\lcross kS_N$.

\begin{corol} The full bosonisation Hopf algebra $\CE_N\lbiprod
(k(S_N)\lcross kS_N)$ contains $\CE_N\lbiprod kS_N$ in
Corollary~6.3 and $k(S_N)$ as sub-Hopf algebras with  additional
relations
\[ f[ij]=[ij]\, L_{(ij)}(f),\quad gf=\ad_{g^{-1}}(f)g,\quad \forall g\in
S_N,\ f\in k(S_N)\] where $L_g(f)=f(g(\ ))$ and
$\ad_g=L_gR_{g^{-1}}$. The same algebra has another `conjugate'
coproduct containing $\CE_N\lbiprod k(S_N)$ and $kS_N$ as sub-Hopf
algebras.
\end{corol}
\proof The right action of $k(S_N)$ on $\CE_N$ given by the
grading can also be used as a left action. This action and the
action of $kS_N$ is the action of the Drinfeld double
corresponding to the crossed module. We make the semidirect
product by this. The quasitriangular structure $\CR=\sum_{g\in
S_N}\delta_g\tens g$ defines a left coaction
$\Delta_L(e_a)=\CR_{21}.e_a=\sum_{g\in S_N}g\tens\delta_g.e_a=
a\tens e_a$ induced from the action, so the same cross coproduct
as for $\CE_N\lbiprod kS_N$. On the other hand every
quasitriangular Hopf algebra has a conjugate quasitriangular
structure $\bar {\CR}=\CR_{21}^{-1}$. We regard the same algebra
$\CE_N$ developed as a braided group in this opposite braided
category (the opposite braided coproduct looks the same on the
generators $E_N$.) Using $\bar{\CR}$ gives a second induced coaction
$\bar\Delta_L(e_a)=\sum_{g\in S_N} \delta_{g^{-1}}\tens g.e_a$.
This gives a second ordinary coproduct
\[\bar\Delta [ij]=[ij]\tens 1+\sum_{g\in
S_N}\delta_{g^{-1}}\tens[g(i),g(j)]\] which is a left handed
version $\CE_N\lbiprod k(S_N)$ of the second biproduct
bosonisation in Corollary~6.3. This is an example of a general
theory in \cite{Ma:qsta} where the two coproducts are related by
complex conjugation in a $*$-algebra setting over $\C$.
 \eproof

Having understood the structure of $\CE_N$ in a natural way, let
us note now that all of the above applies equally well to the full
quotient of it \eqn{Ew}{\CE_w=T E_N/\oplus_n \ker {Sym}_n,\quad
{Sym}_n=\sum_{\sigma=s_{i_1}\cdots s_{i_{l(\sigma)}}\in
S_n}\Psi_{i_1}\cdots\Psi_{i_{l(\sigma)}}} where in principle there
could be nonquadratic relations. In this case, since
${Sym}_n=[n,\Psi]!$, it is clear that here the pairings are now
nondegenerate (we have divided by the coradicals of the pairing in
Proposition~6.4). Therefore $\CE_w$ is a self-dual braided group.
If finite-dimensional then it would inherit a symmetric Hilbert
series as explained in Section~3. Also for the reasons given
there, we expect that $\CE_N$ and $\CE_w$ coincide and the latter
if finite dimensional will have a symmetric Hilbert series which
will prove the conjecture of a symmetric Hilbert series for
$\CE_N$ made in \cite{FomKir:qua}. But if they do not coincide, we
propose $\CE_w$ as the better-behaved version of $\CE_N$; it may
be that $\CE_w$ is finite-dimensional while the $\CE_N$ is likely
not to be for $N\ge 6$. Thus we propose a potential and better
behaved quotient of $\CE_N$. Moreover, our braided group methods
work for general finite groups where we would not expect $\CE_w$
to be quadratic and which would probably be needed for flag
varieties associated to different Lie algebras beyond $SL_N$. This
is a proposal for further work.

\subsection*{Acknowledgements} I would like to thank
S. Fomin and A. Zelevinsky for suggesting to compare with the
algebra $\CE_N$ after a presentation of \cite{Ma:rief} at the
Erwin Schroedinger Institute in 2000. I also want to thank R.
Marsh for the comment about the preprojective algebra after a
presentation in Leicester in 2001, and G. Lehrer for the comment
about cohomology of configuration spaces on a recent visit to
Sydney. The work itself was presented at the Trieste/SISSA
conference, March 2001, at the Banach Center quantum groups
conference in September 2001 and in part at the present
conference. The article was originally submitted to J. Pure and
Applied Algebra in August 2001 and archived on math.QA/0105253;
since then the introduction was redone and some of the more
technical material was moved to an Appendix. The author is a Royal
Society University Research Fellow.

\appendix
\section{Braided group structure of $\Lambda_w$}

Here we will say a little more about the general theory behind
exterior algebras $\Lambda_{quad}$ or $\Lambda_w$ than covered in
the Preliminaries in Section~2. This is needed for some of the
remarks about Hodge * operator mentioned in Sections~3,5 and is
also the motivation behind the results given directly for $\CE_N$
in Section~6. It was considered too technical to be put in the
main text.

First of all, the Woronowicz construction $\Omega_w$ on a quantum
group $A$ is usually given as a quotient of the tensor algebra on
$\Omega^1$ over $A$. We have instead moved everything over to the
left-invariant forms $\Lambda_w$ which is a `braided approach' to
the exterior algebra in \cite{Ma:eps}\cite{BesDra:dif}. See also
\cite{Ros}. The starting point is that associated to any linear
space $\Lambda^1$ equipped with a Yang-Baxter or braid operator
(in our case $-\Psi$) one has $\Lambda_{quad}$ (and similarly
$\Lambda_w$) braided linear spaces with additive coproduct
\eqn{LamDelta}{ \und\Delta e_a=e_a\tens 1+1\tens e_a,\quad
\und\eps e_a=0.} In our case these live in the braided category
which is a $\Z_2$ extension of the category of $A$-crossed
modules, with $\Lambda^1$ odd. Thus one may verify:
\[ \und\Delta (e_ae_b)=(e_a\tens 1+1\tens e_a)(e_b\tens 1+1\tens
e_b)=e_ae_b\tens 1+1\tens e_ae_b+(\id-\Psi)(e_a\tens e_b).\]  If
$\lambda_{a,b}e_ae_b=0$ (summation understood) then $\und\Delta$
of it is also zero since the relation in degree 2 is exactly that
$(\id-\Psi)(\lambda_{a,b}e_a\tens e_b)=0$. This covers
$\Lambda_{quad}$. For $\Lambda_w$ one has to similarly look at the
higher degrees. Similarly to Section~6 there is then a
super-biproduct bosonisation theorem which yields $\Omega_{quad}$
and $\Omega_w$ as super-Hopf algebras by crossed module
constructions. We also have super-braided-partial derivatives
$D_a,\bar D_a$, which define interior products\cite{Ma:eps}.

Here we would like to say a little more as an explanation of the
definition of $\Lambda_w$. Let $\Lambda^{1*}$ be the crossed
module with adjoint braiding $\Psi^*$. It has its own algebra of
`skew invariant tensor fields' \eqn{Lam*}{
\Lambda^*_{quad}=T\Lambda^{1*}/\ker(\id-\Psi^*),\quad
\Lambda^*_w=T\Lambda^{1*}/\oplus_n A_n^*.}

\begin{propos} The tensor algebras $T\Lambda$ and $T\Lambda^{1*}$
are dually paired braided groups as induced by the pairing in
degree $1$, and $\Lambda_w,\Lambda_w^{1*}$ are their quotients by
the kernel of the pairing. \end{propos} \proof Let $\{f^a\}$ be
the dual basis of $\Lambda^{1*}$. The pairing between monomials in
the tensor algebra is then
\[ \<f^{a_n}\cdots f^{a_1},e_{b_1}\cdots
e_{b_m}\>=\delta_{n,m}[n,-\Psi]!^{a_1\cdots a_n}_{b_1\cdots b_n}\]
as for any braided linear space \cite{Ma:book}. In view of
Proposition~2.1 we are therefore defining $\Lambda_w$ exactly by
killing the kernel of the pairing from that side. Similarly from
the other side. \eproof

This the meaning of the Woronowicz construction is that one adds
enough relations that the pairing with its similar dual version is
non-degenerate. Moreover, as in Section~6, we know that by the
theory of integrals on braided groups, if $\Lambda_w$ is finite
dimensional then there is a unique top form ${\rm Top}$, of degree
$d$ say. In this case there is an approach to a Hodge * pairing in
\cite{Ma:eps} based on braided-differentiation of the top form and
related to braided Fourier transform. A similar and more explicit
version of this which has been used in \cite{MaRai:ele} to define
an `epsilon tensor' by $e_{a_1} \cdots e_{a_d}=\eps_{a_1\cdots
a_d}{\rm Top}$ and then use this to define a map
$\Lambda_w^{m}\to\Lambda_w^{*(d-m)}$.  In the presence of an
invariant metric we have $\Lambda^1\isom\Lambda^{*1}$ as
 crossed modules and hence isomorphisms of their generated braided
groups. In this case we have a Hodge
* operator $\Lambda_w^m\to\Lambda_w^{d-m}$. Similarly if
$\Lambda_{quad}$ is finite dimensional.

In the case of a finite group $G$ with calculus defined by a
conjugacy class $\CC$, we compute \eqn{Psi*}{ \Psi^*(f^a\tens
f^b)=f^{a^{-1}ba}\tens f^a} where the adjoint is taken with
respect to the pairing on tensor powers (recall that
conventionally this is defined by pairing the inner factors first
and moving outwards, to avoid unnecessary braid crossings). We let
$\Lambda$ denote either $\Lambda_{quad}$ or $\Lambda_w$ (or an
intermediate quotient).

\begin{corol}If $\CC$ is stable under group inversion then
$\Lambda$ is self-dually paired as a braided group. If $\Lambda$
is finite-dimensional with top degree $d$ we have
\eqn{hodgeG}{*(e_{a_1}\cdots e_{a_m})= =d_m^{-1} \eps_{a_1\cdots
a_d}e_{a_d^{-1}}\cdots e_{a_{m+1}^{-1}}} for some normalisations
$d_m$.
\end{corol}
\proof In this case we have an invariant metric
$\eta^{a,b}=\delta_{a,b^{-1}}$ whereby we identify
$f^a=e_{a^{-1}}$. For $\Lambda_w$ in the algebraically closed case
one would typically chose the $d_m$ so that $*^2=\id$. The formula
as in \cite{MaRai:ele} is arranged to be covariant so that if
${\rm Top}$ is invariant under the $k(G)$-action (which implies
that it commutes with functions) then $*$ will extend to a
bimodule map $\Omega^m\to\Omega^{d-m}$. \eproof

Similarly, the exterior algebra $\Omega$ is generated in the
finite group case by $k(G)$ and $\Lambda$ with the cross relations
(\ref{OmegaG}), which is manifestly a cross product $k(G)\rcross
\Lambda$. The super coalgebra explicitly is \eqn{superG}{ \Delta
e_a=\sum_{g\in G} e_{gag^{-1}}\tens \delta_g+1\tens e_a, \quad
\eps e_a=0} and extends the group coordinate Hopf algebra. Here
$\delta_g$ is a delta-function on $G$. Indeed, the $G$-grading
part of the crossed module structure on $\Lambda^1$ extends to all
of $\Lambda$ and defines a right action of $k(G)$ on it (by
evaluating against the total $G$-degree) which is used in the
cross product algebra. Meanwhile, the left $G$-action defines a
right coaction of $k(G)$, \eqn{rcoaG}{ \Delta_R(e_a)=\sum_g
e_{gag^{-1}}\tens \delta_g} which extends as an algebra
homomorphism to $\Lambda$ because $\Psi$ is $\ad$-covariant.
Semidirect coproduct by this defines the coalgebra of
$k(G)\rbiprod\Lambda$. The two fit together to form a super-Hopf
algebra just because the original structure on $\Lambda^1$ was a
crossed module. For example, one may check
\begin{align*} \Delta(e_a e_b)&=(\sum_g e_{gag^{-1}}\tens\delta_g
+1\tens e_a)(\sum_{h}e_{hbh^{-1}}\tens \delta_h+1\tens e_b)\\
&=1\tens e_ae_b+\sum_g
e_{gag^{-1}}e_{gbg^{-1}}\tens\delta_g+e_{gag^{-1}}
\tens\delta_ge_b-e_{gbg^{-1}}\tens e_a\delta_g\\ &=\left(1\tens
\cdot+\Delta_R\circ\cdot+(\Delta_R\tens
\id)(\id-\Psi)\right)e_a\tens e_b\end{align*} since we are
extending as a super-Hopf algebra (so $\Lambda^1$ is odd). We used
the relations in the algebra and changed variables in the last
term. From this it is clear that $\Delta$ is well defined in the
quotient by $\ker(\id-\Psi)$. This covers $\Lambda=\Lambda_{quad}$
but the same holds also for $\Lambda_w$.

Similarly, since a right $k(G)$-crossed module is the same thing
as a left $kG$-crossed module, we can make another super-Hopf
algebra $\Lambda\lbiprod kG$. We extend the $G$-action
$g.e_a=e_{gag^{-1}}$ to $\Lambda$ for the cross product and the
grading defines a left coaction \eqn{lcoaG}{ \Delta_L e_a=a\tens
e_a} which we extend to products (expressing the total
$G$-degree). Semidirect product and coproduct by these gives
\eqn{dualsupG}{ g e_a=e_{gag^{-1}} g,\quad \Delta e_a=e_a\tens 1 +
a\tens e_a, \quad \eps e_a=0} extending the Hopf algebra structure
of the group algebra $kG$. This time \[ \Delta(e_ae_b)=e_a
e_b\tens 1+ab\tens e_a e_b+e_a b\tens e_b-a e_b\tens e_a\]
\[{\ }\qquad\quad =\left(\cdot\tens
1)+\Delta_L\circ\cdot+(\id\tens\Delta_L)(\id-\Psi)\right)e_a\tens
e_b\] which is well-defined on the quotient . Geometrically, this
is the dual of the super-Hopf algebra $k(G)\rbiprod\Lambda^*$ of
skew-vector fields. These are the direct contructions of the cross
products analogous to those in Section~6 for $\CE_N$. We have
similar pairing results.

Finally, let us note that at this level of generality all the same
proofs work with $-\Psi$ replaced by $\Psi$. Thus for any crossed
module $E$ with braiding $\Psi$ we have a braided space
\eqn{Equad}{ \CE_{quad}=T E/\ker(\id+\Psi)} and similarly $\CE_w$
defined by ${Sym}_n$ as in (\ref{Ew}), both forming additive
braided groups. Moreover, one should be able to construct a
suitable crossed module from any conjugacy class on a finite group
and possibly a cocycle $\zeta$. This indicates how the analogues
of the Fomin-Kirillov algebra could be extended to other types.

\baselineskip 14pt
%\bibliographystyle{unsrt}
%\bibliography{biblio}

\end{document}